\documentstyle[12pt]{amsart}
\expandafter\font\csname cmex/m/n/12\endcsname=cmex10 at 12pt

\expandafter\ifx\csname>krwmacs\endcsname\relax
\expandafter\def\csname>krwmacs\endcsname{done}
\else\endinput\fi
\newtheorem{thm}{Theorem}
\newtheorem{cor}[thm]{Corollary}
\newtheorem{lem}[thm]{Lemma}
\newtheorem{prop}[thm]{Proposition}

\theoremstyle{definition}

\theoremstyle{remark}
\newtheorem{remark}[thm]{Remark}

\newcommand{\thmref}[1]{Theorem~\ref{#1}}

\newcommand{\lemref}[1]{Lemma~\ref{#1}}
\newcommand{\propref}[1]{Prop\-o\-si\-tion~\ref{#1}}
\newcommand{\corref}[1]{Cor\-ol\-lary~\ref{#1}}

\expandafter\ifx\csname >amsstdmathmacs\endcsname\relax
\expandafter\def\csname >amsstdmathmacs\endcsname{done}
\else\endinput\fi
\def\mathcs{C^{\displaystyle *}}
\def\cs{\ifmmode\mathcs\else$\mathcs$\fi}

\def\acg{{A \rtimes_\alpha G}}
\def\Ind{\operatorname{Ind}}

\def\Aut{\operatorname{Aut}}

\def\id{\operatorname{id}}
\def\set#1{\{\,#1\,\}}
\let\tensor=\otimes
\def\({\bigl(}
\def\){\bigr)}

\newbox\hidebox
\def\spechide#1{\setbox\hidebox=\hbox{$#1$}
\hbox to\wd\hidebox{$\box\hidebox^\wedge$\hss}}
\expandafter\ifx\csname>amsspfonts\endcsname\relax
\expandafter\def\csname>amsspfonts\endcsname{done}
\else\endinput\fi
\makeatletter

\new@fontshape{eus}{m}{n}{%
   <5>eusm5%
   <6>eusm6%
   <7>eusm7%
   <8>eusm8%
   <9>eusm9%
   <10>eusm10%
   <11>eusm10 at10.95pt%
   <12>eusm10 at12pt%
   <14>eusm10 at14.4pt%
   <17>eusm10 at17.28pt%
   <20>eusm10 at20.74pt%
   <25>eusm10 at24.88pt}{}

\extra@def{eus}{\skewchar#1'60}{}

\newmathalphabet*{\script}{eus}{m}{n}
\makeatother
\newsymbol\rtimes 226F
\def\bimodfont#1{\script{#1}}
\def\X{{\bimodfont X}}
\def\Y{{\bimodfont Y}}
\def\bZ{{\bimodfont Z}}
\def\gmx{G\backslash X}
\let\xmg=\gmx
\def\xmh{X/H}
\def\hmx{H\backslash X}
\def\Iso{\operatorname{Iso}}
\def\br{\operatorname{Br}}
\def\brg{\br_G}
\def\brh{\br_H}
\def\sbr{{\frak B \frak r}}
\def\sbrg{\sbr_G}
\def\sbrh{\sbr_H}
\def\gmxtensor{\tensor_{C(\gmx)}}

\def\gsigma{\sigma^G}

\def\ttheta{\Theta}
\def\teta{\Lambda}
\def\smeover #1{\,\mathord{\mathop{\text{--}}\nolimits_{#1}}\,}
\def\sme{\,\mathord{\mathop{\text{--}}\nolimits_{\relax}}\,}
\def\eb{im\-prim\-i\-tiv\-ity bi\-mod\-u\-le}
\def\tvo{\tilde v^o}
\def\tv{\tilde v}
\def\ipcomma{\mathrel{,}}
\let\ipscriptstyle=\scriptscriptstyle
\def\lipsqueeze{{\mskip -3.0mu}}
\def\ripsqueeze{{\mskip -3.0mu}}
\newbox\ipstrutbox
\setbox\ipstrutbox=\hbox{\vrule height8.5pt
width 0pt}
\def\ipstrut{\copy\ipstrutbox}
\def\lip#1<#2,#3>{\mathopen{\relax_{\ipstrut\ipscriptstyle{
#1}}\lipsqueeze
\langle} #2\ipcomma #3 \rangle}
\def\blip#1<#2,#3>{\mathopen{\relax_{\ipstrut
\ipscriptstyle{ #1}}\lipsqueeze\bigl\langle} #2\ipcomma #3 \bigr\rangle}
\def\rip#1<#2,#3>{\langle #2\ipcomma #3
\rangle_{\ripsqueeze\ipstrut\ipscriptstyle{#1}}}
\def\brip#1<#2,#3>{\bigl\langle #2\ipcomma #3
\bigr\rangle_{\ripsqueeze\ipstrut\ipscriptstyle{#1}}}
\def\a{\frak a}
\def\b{\frak b}
\def\m(#1,#2){M_{(#1,#2)}}
\def\xm(#1,#2){\X^{\m({#1},{#2})}}
\def\n(#1,#2){N_{(#1,#2)}}
\def\iX{\relax^I\X}
\def\uo{u^o}
\theoremstyle{remark}

\newtheorem{notation}{Notation}

\begin{document}

\title[Isomorphism of Brauer Groups]{The equivariant Brauer
groups of commuting free and proper actions are isomorphic}

\author[Kumjian]{Alexander Kumjian}
\address{Department of Mathematics \\
University of Nevada \\
Reno, NV 89557 \\
USA}
\email{alex@@math.unr.edu}

\author[Raeburn]{Iain Raeburn}
\address{Department of Mathematics \\
University of Newcastle \\
Newcastle, NSW 2308 \\
Australia}
\email{iain@@math.newcastle.edu.au}

\author[Williams]{\\Dana P. Williams}
\address{Department of Mathematics \\
Dartmouth College \\
Hanover, NH 03755-3551 \\
USA}
\email{dana.williams@@dartmouth.edu}
\thanks{The third author was partially supported by the
National Science Foundation.}
\thanks{This research was supported by the
Australian Department of Industry,
Science, and Technology.}



\date{1 September 1994}

\maketitle
%
%
\begin{abstract}
If $X$ is a locally compact
space which admits commuting free and proper actions of
locally compact groups $G$ and $H$,
then the Brauer groups $\brh(\gmx)$ and $\brg(\xmh)$ are
naturally isomorphic.
\end{abstract}
%
%

Rieffel's formulation of Mackey's Imprimitivity Theorem asserts that
if $H$ is a closed subgroup of a locally compact group $G$, then the
group \cs-algebra $\cs(H)$ is Morita equivalent to the crossed product
$C_0(G/H)\rtimes G$. Subsequently, Rieffel found a symmetric version,
involving two subgroups of $G$, and Green proved the following {\it
Symmetric Imprimitivity Theorem\/}: if two locally compact groups act
freely and properly on a locally compact space $X$, $G$ on the left
and $H$ on the right, then the crossed products $C_0(\gmx)\rtimes H$
and $C_0(\xmh)\rtimes G$ are Morita equivalent. (For a discussion and
proofs of these results, see \cite{rieff5}.) Here we shall show
that in this situation there is an isomorphism
$\brh(\gmx)\cong\brg(\xmh)$ of the equivariant Brauer groups
introduced in \cite{ckrw}.

Suppose $(G,X)$ is a second countable locally compact transformation
group. The objects in the underlying set $\sbrg(X)$ of the equivariant
Brauer group $\brg(X)$ are dynamical systems $(A,G,\alpha)$, in which
$A$ is a separable continuous-trace \cs-algebra with spectrum $X$,
and $\alpha:G\to \Aut(A)$ is a strongly continuous action of $G$ on
$A$ inducing the given action of $G$ on $X$. The equivalence relation
on such systems is the equivariant Morita equivalence studied in
\cite{combes}, \cite{cmw}. The group operation is given by
$[A,\alpha]\cdot[B,\beta]=[A\otimes_{C(X)}B,\alpha\otimes\beta]$, the
inverse of $[A,\alpha]$ is the conjugate system $[\overline
A,\overline\alpha]$, and the identity is represented by
$(C_0(X),\tau)$, where $\tau_s(f)(x)=f(s^{-1}\cdot x)$.

\begin{notation}
Suppose that $H$ is a locally compact group, that $X$ is a free
and proper right
$H$-space, and that $(B,H,\beta)$ a dynamical system.  Then
$\Ind_H^X(B,\beta)$ will be the \cs-algebra (denoted by
$GC(X,B)^{\alpha}$ in $\cite{rw}$ and by $\Ind(B;X,H,\beta)$ in
\cite{ra2})
of bounded continuous functions
$f:X\to B$ such that
$\beta_h\(f(x\cdot h)\)=f(x)$, and
$x\cdot H\mapsto\|f(x)\|$ belongs to $C_0(\xmh)$.
\end{notation}
We now
state our main theorem.

\begin{thm}\label{thm:mt}
Let $X$ be a second countable locally compact Hausdorff space, and
let $G$ and $H$ be second countable locally compact groups.
Suppose that $X$ admits a free and proper left $G$-action, and a
free and proper right $H$-action such that $(g\cdot x)\cdot
h=g\cdot(x\cdot h)$ for all $x\in X$, $g\in G$, and $h\in H$.
Then there is an isomorphism $\Theta$ of $\brh(\gmx)$ onto
$\brg(\xmh)$ satisfying:
\begin{enumerate}
\item
  if $(A,\alpha)$ represents $\Theta[B, \beta]$,
then
$A\rtimes_\alpha G$ is Morita equivalent to $B\rtimes_\beta H$;
\item
 $\Theta[B,\beta]$ is
realised by the pair $(\Ind_H^X(B,\beta)/J,\tau\otimes \id)$ in
$\sbrg(\xmh)$, where
$\tau\otimes\id$ denotes left translation
and, if $\pi_{G\cdot x}$ is the element of $\widehat B=\gmx$
corresponding to $G\cdot x$,
\begin{equation*}
J=\set{f\in \Ind_H^X(B,\beta):\text{$\pi_{G\cdot x}\(f(x)\)=0$ for all
$x\in X$}}.
\end{equation*}
\end{enumerate}
\end{thm}

Item~(1) is itself a generalization of Green's symmetric
imprimitivity theorem, and our proof of \thmref{thm:mt} follows
the approach to Green's theorem taken in
\cite{cmw}:
 prove that both $C_0(\gmx)\rtimes H$ and
$C_0(\xmh)\rtimes G$ are Morita equivalent to
$C_0(X)\rtimes_\alpha(G\times H)$, where
$\alpha_{s,h}(f)(x)=f(s^{-1}\cdot x\cdot h)$, by noting that the
Morita equivalences of $C_0(X)\rtimes G$ with $C_0(\gmx)$ and
$C_0(X)\rtimes H$ with $C_0(\xmh)$ (\cite{green2},
\cite[Situation~10]{rieff5}) are equivariant,
and hence induce Morita equivalences
\begin{align*}
C_0(\gmx)\rtimes H & \sim \(C_0(X)\rtimes G\)\rtimes H
\cong C_0(X)\rtimes(G\times H) \\
&\cong \(C_0(X)\rtimes H\)\rtimes G
\sim C_0(\xmh)\rtimes G.
\end{align*}
The same symmetry considerations show that it will be enough to prove that
$\brh(\gmx)\cong\br_{G\times H}(X)$. Since we already know that
$\br(\gmx)\cong\brg(X)$ \cite[\S6.2]{ckrw}, we just have to
check that this isomorphism is compatible with the actions of
$H$.

Suppose $G$ acts freely and properly on $X$, and
$p:X\to\gmx$ is the orbit map.
If $B$ is a
a \cs-algebra with a nondegenerate action of $C_0(\gmx)$, then the
pull-back
$p^*B$ is the quotient of $C_0(X)\otimes B$ by the balancing ideal
\begin{equation*}
I_{\gmx}=\overline{\operatorname{span}}
\{f\cdot \phi\otimes b-\phi\otimes f\cdot b:
\phi\in C_0(X),f\in C_0(\gmx),b\in B\}
\end{equation*}
in other words, $p^*B=C_0(X)\otimes_{C(\gmx)}B$. The nondegenerate
action of $C_0(\gmx)$ on $B$ induces a continuous map $q$ of $\widehat
B$ onto $\gmx$, characterized by $\pi(f\cdot b)=f(q(\pi))\pi(b)$. Then
under the natural identification of $C_0(X)\otimes B$ with $C_0(X,B)$,
\begin{equation*}
I_{\gmx}\cong\{f\in C_0(X,B): \pi(f(x))=0 \hbox{ for all
}x\in q(\pi)\},
\end{equation*}
so that $p^*B$ has spectrum
\begin{equation*}
\widehat{p^*B}=\{(x,\pi)\in X\times\widehat B: G\cdot x=q(\pi)\}.
\end{equation*}
If $B$ is a continuous-trace algebra with spectrum $\gmx$, then $p^*B$ is
a continuous-trace algebra with spectrum $X$.

The isomorphism $\Theta:\br(\gmx)\cong\brg(X)$ is given by
$\Theta[A]=[p^*A,\tau\otimes{\id}]$. To prove $\Theta$ is
surjective in \cite{ckrw}, we used \cite[Theorem 1.1]{rr}, which
implies that if $(B,\beta)\in \sbrg(X)$, then $B\rtimes_\beta G$ is a
continuous-trace algebra with spectrum $\gmx$ such that $(B,\beta)$ is
Morita equivalent to $\(p^*(B\rtimes_\beta G),\tau\otimes{\id}\)$,
and hence that $[B,\beta]=\Theta[B\rtimes_\beta G,{\id}]$.
In obtaining the required
equivariant version of
\cite[Theorem 1.1]{rr}, we have both simplified
the proof and mildly strengthened the conclusion (see \corref{cor:rr}
below). However, with all these different group actions around, the
notation could get messy, and we pause to establish some conventions.

\begin{notation}
We shall be dealing with several spaces carrying a left action of $G$
and/or a right action of $H$. We denote by $\tau$ the action of $G$ by
left translation on $C_0(G)$, $C_0(X)$ or $C_0(\gmx)$, and by $\sigma$
any action of $H$ by right translation; we shall also use $\sigma^G$
to denote the action of $G$ by right translation  on $C_0(G)$.
Restricting an action $\beta$ of $G\times H$ on an algebra $A$ gives
actions $\alpha:G\to \Aut(A)$, $\gamma:H\to\Aut(A)$ such that
\begin{equation}\label{eq:commute}
\alpha_s\(\gamma_h(a)\)=\gamma_h\(\alpha_s(a)\)
\quad\text{for all $h\in H$, $s\in G$, $a\in A$.}
\end{equation}
Conversely, two actions $\alpha,\gamma$ satisfying
\eqref{eq:commute} define an action of $G\times H$ on $A$, which we
denote by $\alpha\gamma$; we write $\gamma$ for $\id\gamma$ since
it will be clear from context whether an action of $H$ or or
$G\times H$ is called for. If
$\Phi:(A,G,\alpha)\to(B,G,\beta)$ is an equivariant isomorphism (i.e.
$\Phi(\alpha_s(a))=\beta_s(\Phi(a))\,$), then we denote by
$\Phi\rtimes\id$ the induced isomorphism of $A\rtimes_\alpha G$ onto
$B\rtimes_\beta G$. Similarly, if
$\alpha$ and $\gamma$ satisfy
\eqref{eq:commute}, we write $\alpha\rtimes\id$ for the induced action
of $G$ on $A\rtimes_\gamma H$.
\end{notation}

\begin{lem}\label{lem:1}
Suppose a locally compact group $G$ acts freely and properly
on a locally compact space $X$, and that $A$ is a \cs-algebra carrying
a non-degenerate action of
$C_0(X)$. If $\alpha:G\to \Aut(A)$ is an action of $G$ on $A$
satisfying $\alpha_s(\phi\cdot a)=\tau_s(\phi)\cdot\alpha_s(a)$, then
the map sending $f\tensor a$ in $C_0(X)\tensor A$ to the function
$s\mapsto f\cdot\alpha_s^{-1}(a)$ induces an equivariant isomorphism
$\Phi$ of $\(C_0(X)\gmxtensor A,G,\id\otimes\alpha\)$ onto
$\(C_0(G,A),G,\tau\otimes\id\)$.
\end{lem}

\begin{remark}\label{rmk} For motivation, consider the case where
$A=C_0(X)$. Then the map $\Psi:C_b(X\times X)\to C_b(G\times X)$
defined by
$\Psi(f)(s,x)=f(x,s\cdot x)$ maps $C_0$ to $C_0$ precisely when the
action is proper, has range which separates the points of $G\times X$
precisely when the action is free, and has kernel consisting of the
functions which vanish on the closed subset $\Delta=\{(x,y):G\cdot
x=G\cdot y\}$. Thus the free and proper actions are precisely those
for which $\Psi$ induces an isomorphism of
$C_0(X)\otimes_{C(\gmx)}C_0(X)$ onto $C_0(G)\otimes C_0(X)$.
\end{remark}

\begin{pf*}{Proof of Lemma \ref{lem:1}}
If $\phi\in C_0(\gmx)$, then
$f\cdot\phi
\tensor a$ and $f\tensor \phi\cdot a$ have the same image in
$C_0(G,A)$, and the map factors through the balanced tensor product as
claimed. Further, $\Phi$ is related to the map $\Psi$ in Remark~\ref{rmk}
by
\begin{equation}\label{eq:phi}
\Phi(f\tensor g\cdot a) =
\(\Psi(f\otimes g)(s,\cdot)\)\cdot\alpha_s^{-1}(a).
\end{equation}
Thus it follows from the remark that \eqref{eq:phi} defines an element
of
$C_0(G,A)$ and that the closure of the range of $\Phi$ contains all
functions of the form $s\mapsto \xi(s)f\cdot\alpha_s^{-1}(a)$
for $\xi\in C_c(G)$, $f\in C_c(X)$, and $a\in A$.  These
elements span a dense subset of $C_0(G,A)$, and hence $\Phi$ is
surjective. The nondegenerate action of $C_0(X)$ on $A$ induces a
continuous equivariant map $q$ of $\widehat A$ onto $X$ such that
$\pi(f\cdot a)=f(q(\pi))\pi(a)$, and the balanced tensor product
$C_0(X)\otimes_{C(\gmx)}A$ has spectrum
$\Delta=\{(x,\pi):G\cdot x=G\cdot q(\pi)\}$. Since each
representation $\(q(\pi),s\cdot
\pi\)=\(q(\pi),\pi\circ\alpha_s^{-1}\)$ in $\Delta$ factors through
$\Phi$ and the representation $b\mapsto \pi(b(s))$ of $C_0(G,A)$, the
homomorphism $\Phi$ is also injective. Finally, to see the
equivariance, we compute:
\begin{equation*}
\Phi\(\id\otimes \alpha_s(h\otimes
a)\)(t)=h\cdot\alpha_t^{-1}\(\alpha_s(a)\)
=\Phi(h\otimes a)(s^{-1}t)=\tau_s\otimes\id\(\Phi(h\otimes a)\)(t).
\qed
\end{equation*}
\renewcommand{\qed}{}\end{pf*}

\begin{cor}\label{cor:rr} {\normalshape(cf. \cite[Theorem 1.1]{rr})} Let
$(G,X)$ and $\alpha:G\to \Aut(A)$ be as in Lemma \ref{lem:1}. Then
there is an equivariant isomorphism of $\(p^*(A\rtimes_\alpha
G),G,p^*\id\)$ onto $\(A\otimes{\cal
K}(L^2(G)),G,\alpha\otimes{\normalshape{\text {Ad}}}\,\rho\)$.
\end{cor}

\begin{pf}
A routine calculation shows that
the equivariant isomorphism $\Phi$ of Lemma \ref{lem:1} gives an
equivariant isomorphism
\begin{multline}\label{eq:added}
\Phi\rtimes\id:\(\(C_0(X)\otimes_{C(\gmx)}A\)\rtimes_{\id\otimes\alpha}G,
(\tau\otimes\id)\rtimes\id\) \\
\to\(C_0(G,A)\rtimes_{\tau\otimes\id}G,
(\sigma^G\otimes\alpha)\rtimes\id\).
\end{multline}
We also have equivariant isomorphisms
\begin{equation}\label{420}
\begin{split}
\(C_0(G,A)\rtimes_{\tau\otimes\id}G,
(\sigma^G\otimes\alpha)\rtimes\id\)&\cong\(A\otimes\(C_0(G)\rtimes_\tau
G\),\alpha\otimes(\sigma^G\rtimes\id)\),\\
&\cong\(A\otimes{\cal K}(L^2(G)),\alpha\otimes{\text
{Ad}}\,\rho\)\\
\end{split}
\end{equation}
and
\begin{equation}\label{421}
\(C_0(X)\otimes_{C(\gmx)}(A\rtimes_\alpha G),\tau\otimes\id\)\cong
\(\(C_0(X)\otimes_{C(\gmx)}A\)\rtimes_{\id\otimes\alpha}G,
(\tau\otimes\id)\rtimes\id\);
\end{equation}
combining \eqref{eq:added}, \eqref{420}, and \eqref{421}
gives the result.
\end{pf}

\begin{lem}\label{lem:2}
In addition to the hypotheses of
\lemref{lem:1}, suppose that $H$  is a locally
compact group acting on the right of $X$,
and that
$(A,H,\gamma)$ is a dynamical system such that
$\alpha$ and
$\gamma$ commute and
$\gamma_h(f\cdot
a)=\sigma_h(f)\cdot\gamma_h(a)$ for $h\in H$, $f\in C_0(X)$,
$a\in A$.
Then the action $\tau\sigma\tensor\gamma$ of
$G\times H$ on $C_0(X)\tensor
A$ preserves the balancing ideal $I_{\gmx}$, and
hence induces an action of $G\times H$ on
$C_0(X)\gmxtensor A$, also denoted $\tau\sigma\tensor\gamma$.
The equivariant isomorphism of Lemma \ref{lem:1}
induces an equivariant
isomorphism
\begin{equation*}
\(\(C_0(X)\gmxtensor
A\)\rtimes_{\id\tensor\alpha}G,(\tau\sigma\tensor\gamma)\rtimes\id\)
\cong\(C_0(G,A)
\rtimes_{\tau\tensor\id}G,(\gsigma\tensor\alpha
\gamma)\rtimes\id\).
\end{equation*}
\end{lem}
\begin{pf}
The first assertion is straightforward.
For the second,
we can consider the actions of $H$ and $G$ separately.
We have already observed in \eqref{eq:added} that
$\Phi\rtimes\id$ intertwines the $G$-actions.
On the other hand, if $h\in H$ and $t\in G$, then
\begin{align*}
\Phi\(\sigma_h\tensor\gamma_h(f\tensor a)\)(t) =&
\sigma_h(f)\cdot\alpha_t^{-1}\(\gamma_h(a)\)
=\sigma_h(f)\cdot\gamma_h\(\alpha_t^{-1}(a)\) \\
&=\gamma_h\(\Phi(f\tensor a)(t)\).
\qed
\end{align*}
\renewcommand{\qed}{}\end{pf}

\begin{cor}\label{cor:42}
Let ${}_GX_H$ and $\alpha:G\to\Aut(A)$, $\gamma:H\to\Aut(A)$ be as in
the lemma. Denote by $p$ the orbit map of $X$ onto $\gmx$. Then there is
an equivariant isomorphism
\begin{equation*}
\(p^*(A\rtimes_\alpha G),G\times H,\tau\sigma\otimes(\gamma\rtimes\id)\)
\cong\(A\otimes{\cal K}(L^2(G)), G\times
H,\alpha\gamma\otimes{\normalshape {\text{Ad}}}\,\rho\).
\end{equation*}
\end{cor}
\begin{pf}
Compose the isomorphism of Lemma \ref{lem:2} with
(\ref{420}) and (\ref{421}).
\end{pf}

We are now ready to define our map of $\brh(\gmx)$ into $\br_{G\times
H}(X)$. Suppose $(B,\beta)\in\sbrh(X)$. Then the action
$\tau\sigma\otimes\beta$ of $G\times H$ preserves the balancing ideal
$I_{\gmx}$: if $\phi\in C_0(\gmx)$ then
\begin{align*}
(\tau\sigma\otimes\beta)_{s,h}(f\cdot\phi\tensor b&-f\tensor \phi\cdot b)=
\sigma_h\(\tau_s(f\cdot\phi)\)\tensor\beta_h(b) -
\sigma_h\(\tau_s(f)\)\otimes\beta_h(\phi\cdot b) \\
&=\sigma_h\(\tau_s(f)\)\cdot\sigma_h(\phi)\tensor \beta_h(b) -
\sigma_h\(\tau_s(f)\)\otimes\sigma_h(\phi)\cdot\beta_h(b).
\end{align*}
Since $p^*(B)$ is a continuous-trace \cs-algebra with spectrum
$X$ \cite[Lemma~1.2]{rr}, and $\tau\sigma\otimes\beta$
covers the canonical $G\times H$-action on $X$,
 we can define $
\theta:\sbrh(\gmx)\to\sbr_{G\times H}(X)$ by
$\theta(B,\beta)=(p^*(B),\tau\sigma\otimes\beta)$.

Similarly if $(A,\alpha\gamma)\in\sbr_{G\times H}(X)$, then
$\acg$ is a continuous-trace \cs-algebra with spectrum $\gmx$ by
\cite[Theorem~1.1]{rr}. Since $\gamma$ is compatible with $\sigma$, we
have
$\gamma_h\(\phi\cdot z(s)\)=\sigma_h(\phi)\cdot\gamma_h\(z(s)\)$ for
$z\in C_c(G,A)$, and hence $\gamma\rtimes\id$ covers the given action of
$H$ on $X$. Thus we can define $\lambda:\sbr_{G\times
H}(X)\to\sbrh(\gmx)$ by
$\lambda(A,\alpha\gamma)=(\acg,\gamma
\rtimes \id)$.

\begin{prop}\label{prop:pone}
Let $X$ be a second countable locally compact Hausdorff space, and
let $G$ and $H$ be second countable locally compact groups.
Suppose that $X$ admits a free and proper left $G$-action, and an
$H$-action such that $(g\cdot x)\cdot h = g\cdot (x\cdot h)$ for all
$x\in X$, $g\in G$, and $h\in H$. Then
$\theta$ and
$\lambda$ above preserve Morita equivalence classes, and define
homomorphisms
$\ttheta:\brh(\gmx)\to
\br_{G\times H}(X)$ and $\teta:\br_{G\times H}(X)\to\brh(\gmx)$.
In fact, $\ttheta$ is an isomorphism with inverse $\teta$, and if
$\ttheta[B,\beta]=[A,\alpha]$, then $B\rtimes_\beta H$ is Morita
equivalent to $A\rtimes_\alpha (G\times H)$.
\end{prop}
\begin{pf}
If $(\Y,v)$ implements an equivalence between $(B,\beta)$
and $(B',\beta')$ in $\sbrh(\gmx)$,  then,
the external tensor product
$\bZ=C_0(X)\widehat\tensor\Y$, as defined in \cite[\S1.2]{je-th}
or \cite[\S2]{ckrw}, is a $C_0(X)\tensor B\sme C_0(X)\tensor B'$-\eb.
A routine argument, similar to that in \cite[Lemma~2.1]{ckrw},
shows that the Rieffel correspondence \cite[Theorem~3.1]{rieff2}
between the lattices of ideals in $C_0(X)\tensor
B$ and in $C_0(X)\tensor B'$ maps the balancing ideal $I=I_{
C(\gmx)}$ in $C_0(X)\tensor
B$ to the balancing ideal $J=J_{C(\gmx)}$ in $C_0(X)\tensor B'$.
Thus \cite[Corollary~3.2]{rieff2} implies that $\X=\bZ/\bZ\cdot
J$ is a $p^*(B)\sme p^*(B')$-\eb.
Since $f\cdot
x=x\cdot f$ for all $x\in \X$ and $f\in C_0(X)$, it follows from
\cite[Proposition~1.11]{ra1} that $\X$ implements a Morita
equivalence over $X$.  More tedious but routine calculations show that
the map defined on elementary tensors in $\bZ_0=C_0(X)\odot\Y$ by
$u^o_{(s,h)}(f\tensor y)=\sigma_h\(\tau_s(f)\)\tensor v_h(y)$
extends to the completion $\bZ$, and defines a strongly continuous
map $u:G\times H\to\Iso(\X)$ such that $(\X,u)$ implements an
equivalence between $(p^*(B),\tau\sigma\otimes\beta)$ and
$(p^*(B'),\tau\sigma\otimes\beta')$. Thus $\ttheta$ is well-defined.

Observe that
\begin{align}
\ttheta\([B,\beta][B',\beta']\) &= \ttheta\([ B\gmxtensor B',
\beta\tensor\beta']\)\notag\\
&=
[p^*\(B\gmxtensor B'\),\tau\sigma\tensor(\beta\tensor\beta')].
\label{eq:prod}
\end{align}
But
\eqref{eq:prod} is the class of
\begin{align*}
\(C_0&(X)\gmxtensor B\gmxtensor
B',\tau\sigma\tensor\beta\tensor
\beta'\)\\
&\sim
\(C_0(X)\tensor_{C(X)} C_0(X)\gmxtensor B \gmxtensor B',
\tau\sigma\tensor
\tau\sigma\tensor\beta\tensor
\beta'\) \\
&\sim
\(C_0(X)\gmxtensor B \tensor_{C(X)} C_0(X) \gmxtensor B',
\tau\sigma\tensor \beta
\tensor \tau\sigma\tensor \beta'\),
\end{align*}
which represents the product of $\ttheta[B,\beta]$ and
$\ttheta[B',\beta']$. Thus
$\ttheta$ is a homomorphism.

Now suppose that $(A,\alpha\gamma)\sim(A',\alpha'\gamma'
)$ in $\sbr_{G\times H}(X)$ via $(\bZ,w)$.  Then $u_s=w_{(s,e)}$
and $v_h=w_{(e,h)}$ define actions of $G$ and $H$, respectively,
on $\bZ$.  In particular, $(\bZ,u)$ implements an equivalence
between $(A,\alpha)$ and $(A',\alpha')$ in $\sbr_G(X)$.
It follows from \cite[\S6]{combes} that $\Y_0=C_c(G,\bZ)$
can be completed to a $\acg\sme A'\rtimes_{\alpha'}G$-\eb{} $\Y$.
One can verify that the induced $C_0(\gmx)$-actions on $\Y_0$ are
given by $(\phi\cdot x)(t)=\phi\cdot\(x(t)\)$ and $(x\cdot
\phi)(t) = \(x(t)\)\cdot\phi$, and \cite[Proposition 1.11]{ra1} implies
that
$\Y$ is an
\eb{} over
$\gmx$.  Now define $\tvo_h$ on $\Y_0$ by $\tvo_h(x)(t)=
v_h\(x(t)\)$.  Using the inner products defined in \cite[\S6]{combes},
\begin{align*}
\blip \acg <\tvo_h(x),\tvo_h(y)>(t)
&= \int_G\blip
A<\tvo_h(x)(s),\Delta(t^{-1}s)u_t\(\tvo_h(y)(t^{-1}s)\)>\,ds \\
&=\int_G\blip A <v_h\(x(s)\), \Delta(t^{-1}s)u_t\( v_h\(y(t^{-1}
s)\)\)>\, ds\\
&= \gamma_h\(\lip\acg<x,y>(t)\),
\end{align*}
where, in the last equality, we use $u_s\circ
v_h=v_h \circ u_s$.
A similar computation shows that $\brip A'\rtimes_{\alpha'}G
<\tvo_h(x),\tvo_h(y)>(t)=\gamma'_h\(\rip A'\rtimes_{\alpha'}G
<x,y>(t)\)$.  Thus $\tvo_h$ extends to all of $\Y$ and defines a
map $\tv:H\to\Iso(\Y)$, and it is not hard to verify that $\tv$ is
strongly continuous.  Therefore $(\acg,\gamma\rtimes\id)\sim
(A'\rtimes_{\alpha'}G,\gamma'\rtimes\id)$ in $\sbr_{G\times H}(X)$,
and $\teta$ is well-defined.

Now it will suffice to show that, for $\a\in\sbrh(\gmx)$ and $\b\in
\sbr_{G\times H}(X)$, $\theta\(\lambda(\b)\)\sim\b$ and $\lambda\(
\theta(\a)\)\sim \a$.  For the first of these,
suppose that
$(A,\alpha\gamma)\in\sbr_{G\times H}(X)$.
Then
$\theta\(\lambda(A,\alpha\gamma)\)=\(p^*(\acg),
(\tau\sigma\tensor\gamma)\rtimes\id\)$, which by Corollary \ref{cor:42}
is equivalent to
$\(A\otimes {\cal K}(L^2(G)),\alpha\gamma\otimes{\text {Ad}}\,\rho\)$, and
hence to $(A,\alpha\gamma\)$. For the other direction, suppose that
$(B,\beta)\in\sbrh(\gmx)$. Then
$\lambda\(\theta(B,\beta)\)=\(p^*B\rtimes_{\tau\otimes\id}G,
(\sigma\tensor\beta)\rtimes\id\)$. Now
\begin{equation*}
p^*B\rtimes_{\tau\otimes\id}G\cong
\(C_0(X)\otimes_{C(\gmx)}B\)\rtimes_{\tau\otimes\id}G
\cong\(C_0(X)\rtimes_\tau G\)\otimes_{C(\gmx)}B,
\end{equation*}
which is Morita equivalent to $C_0(\gmx)\otimes_{C(\gmx)}B\cong B$.
Because the Morita equivalence of
$C_0(X)\rtimes G$ with $C_0(\gmx)$ is $H$-equivariant \cite{cmw},
it follows
that
\begin{equation*}
\lambda\(\theta(B,\beta)\)=\(p^*B\rtimes_{\tau\otimes\id}G,
(\sigma\tensor\beta)\rtimes\id\)\sim
\(C_0(\gmx)\otimes_{C(\gmx)}B,\sigma\otimes\beta\)\cong(B,\beta).
\end{equation*}
This shows that $\Lambda\circ\Theta$ is the identity, and also implies
that
\begin{equation*}
p^*B\rtimes_{\tau\sigma\otimes\beta}(G\times
H)\cong\(p^*B\rtimes_{\tau\otimes\id}G\)\rtimes_{\sigma\otimes\beta}H
\cong B\rtimes_\beta H,
\end{equation*}
which proves the last assertion.
\end{pf}

\begin{remark}
We showed that $\Lambda$ is a well-defined map of $\br_{G\times H}(X)$
into $\brh(\gmx)$, and that it is a set-theoretic inverse for $\Theta$;
since $\Theta$ is a group homomorphism, it follows that $\Lambda$ is
also a homomorphism. This seems to be non-trivial:
it implies that if $(A,\alpha)$, $(B,\beta)$ are
in $\sbrg(X)$, then
$(A\otimes_{C(X)}B)\rtimes_{\alpha\otimes\beta}G$ is Morita equivalent to
$(A\rtimes_\alpha G)\otimes_{C(\gmx)}(B\rtimes_\beta
G)$.
We do not know what general mechanism is at work
here. Certainly, it is a Morita equivalence rather than an isomorphism:
if $G$ is finite and the algebra commutative, one algebra is
$|G|$-homogeneous and the other
$|G|^2$-homogeneous. The only
direct way we have found uses
\cite[Theorem 17]{green1}, which seems an excessively heavy
sledgehammer.
\end{remark}

\begin{pf*}{Proof of \thmref{thm:mt}}
It follows from \propref{prop:pone} that there are isomorphisms
$\Theta_H:\brh(\gmx)\to\br_{G\times H}(X)$, and
$\Lambda_G:\br_{G\times H}(X)\to\brg(\xmh)$.
Therefore $\Lambda_G\circ\Theta_H$ is an
isomorphism of $\brh(\gmx)$ onto $\brg(\xmh)$.
Assertion (1)
also follows from Proposition \ref{prop:pone}. The isomorphism
$\Lambda_G\circ\Theta_H$ maps the class of
$(B,\beta)$ in
$\sbrh(\xmg)$ to the class of
$\(p^*(B)\rtimes_{\sigma\tensor\beta}H,(\tau\tensor\id)
\rtimes\id\)$, so it remains to show that the latter
is equivalent to $(A/J,\tau)$.

For convenience, write $I$ for the balancing ideal
$I_{C(\gmx)}$ in $C_0(X)\tensor B$.  Then
\begin{equation*}
p^*(B)\rtimes_{\sigma\tensor\beta}H = \(\(C_0(X)\tensor B\)/I\)
\rtimes_{\sigma\tensor\beta}H =
\(C_0(X,B)\rtimes_{\sigma\tensor\beta}H\) /
\(I\rtimes_{\sigma\tensor\beta}H\)
\end{equation*}
by, for example, \cite[Proposition~12]{green1}.
By \cite[Theorem~2.2]{rw}, $\X_0=C_c(X,B)$
can be completed to a
$C_0(X,B)\rtimes_{\sigma\tensor\beta}H \sme A$-\eb{} $\X$.
The irreducible representations of $A$ are given by
$\m(x,\pi_{G\cdot y})(f)(x)=\pi_{G\cdot y}\(f(x)\)$
\cite[Lemma~2.6]{rw}.
In the proof of \cite[Theorem~2.5]{rw}, it was
shown that the representation $\xm(x,\pi_{G\cdot y})$ of $
C_0(X,B)\rtimes_{\sigma\tensor\beta}H$ induced from $\m(x,\pi_{G\cdot
y})$ via $\X$ is equivalent to $\Ind_{\set e}^G\n(x,G\cdot y)$,
where $\n(x,G\cdot y)$ is the analogous irreducible representation
of $C_0(X,B)$.
Since the orbit space for a proper action is Hausdorff,
\cite{effros} implies that
$\(C_0(X,B),H,{\sigma\tensor\beta}\)$ is regular.
Since $R=\bigoplus_{x\in X}\n(x,G\cdot x)$ is a faithful
representation of $p^*(B)$, it follows from
\cite[Theorem~24]{green1} that $\Ind_{\set e}^G(R)$ is a faithful
representation of $p^*(B)\rtimes_{\sigma\tensor\beta}H$, and so
has kernel $I\rtimes_{\sigma\tensor\beta}H$.
On the other hand, $\Ind_{\set e}^G(R)$ is equivalent
to $\bigoplus_{x\in X} \xm(x,G\cdot x)$.  It follows from
\cite[\S3]{rieff2} that $\iX=\X/I\cdot\X$ is an $p^*(B)
\rtimes_{\sigma\tensor\beta}H\smeover {\xmh}A/J$-\eb.
Then the map $\uo_s:\X_0\to\X_0$ defined by
$\uo_s(\xi)(x)=\xi(s^{-1}
\cdot x)$ induces a map $u:G\to\Iso(\iX)$ such that $(\iX,u)$
implements the desired equivalence.
\end{pf*}

We close with two interesting special cases where the isomorphism
takes a particularly elegant form.
Recall that
if $B$ is a continuous-trace \cs-algebra with spectrum $X$,
then we may view $B$ as the sections
$\Gamma_0(\xi)$ of a \cs-bundle $\xi$ vanishing at infinity.

\begin{cor}
Suppose that $H$ is a closed subgroup of a second countable locally
compact group $G$, and that $X$ is a second countable locally compact
right $H$-space.  Then $G\times X$ is a free and proper $H$-space via the
diagonal action $(s,x)\cdot h=(sh,x\cdot h)$.  Thus $(G\times X)/H$ is a
locally compact $G$-space via $s\cdot[r,x]=[sr,x]$, and the map
$(B,\beta)\mapsto\(\Ind_H^G(B,\beta),\tau\)$ induces an isomorphism of
$\br_H(X)$ onto $\br_G\((X\times G)/H\)$.
\end{cor}
\begin{pf}
We apply  \thmref{thm:mt} to ${}_G(G\times X)_H$, where $G$ acts
on the left of the first factor, obtaining an isomorphism of
$\br_H(X)\cong \br_H\(G\backslash (G\times
X)\)$ onto $\br_G\((G\times X)/H\)$ sending the class of $(B,\beta)$ to
the class of $\Ind_H^{G\times X}(B,\beta)/J$ where
$J=\set{f:f(s,x)(x)=0}$.\looseness=-1

Given $f\in\Ind_H^{G\times X}(B,\beta)$ and $s\in G$, let $\Phi(f)(s)$
be the function from $X$ to $\xi$ defined by $\Phi(f)(s)(x)=f(s,x)(x)$.
We claim $\Phi(f)(s)\in\Gamma_0(\xi)$.  If $x_0\in X$, then $x\mapsto
f(s,x_0)(x)$ is in $\Gamma_0(\xi)$, and $\|\Phi(s)(x)-f(s,x_0)(x)\|$
tends to zero as $x\to x_0$.  It follows from \cite[Proposition~1.6
(Corollary~1)]{fell4} that
$\Phi(f)(s)$ is continuous.  To see that $\Phi(f)(s)$ vanishes at
infinity, suppose that $\set{x_n}\subset X$ satisfies
$$\|\Phi(f)(s)(x_n)\|\ge\epsilon>0$$
for all $n$.  Then $\|f(s,x_n)\|\ge\epsilon$ for all $n$, and passing to
a subsequence and relabeling if necessary, there must be $h_n\in H$ such
that $(s\cdot h_n,x_n\cdot h_n)\to (r,x)$.  Then $h_n\to s^{-1}r\in H$,
and $x_n\to x\cdot(r^{-1}s)$.  In sum, $\Phi(f)(s)\in\Gamma_0(\xi)=B$.
Now the continuity of $f$ easily implies that $s\mapsto\Phi(f)(s)$ is
continuous from $G$ to $B$.
Furthermore, since $\beta$ covers $\sigma$ (i.e., $\beta_h(\phi\cdot
b)(x)=\phi(x\cdot h)\beta_h(b)(x)$)
,
$$f(rh,x)(x)=\beta_h^{-1}\(\Phi(f)(r)\)(x),$$
and $\Phi$ is a $*$-homomorphism of $\Ind_H^{G\times X}(B,\beta)$ into
$\Ind_H^G(B,\beta)$, which clearly has kernel $J$.

Finally, it is not difficult
(cf., e.g., \cite[Lemma~2.6]{rw})
to see that $\Phi\(\Ind_H^{G\times X}(B,\beta)\)$ is a rich subalgebra
of $\Ind_H^G(B,\beta)$ as defined in \cite[Definition~11.1.1]{dix}.
Thus $\Phi$ is surjective by \cite[Lemma~11.1.4]{dix}.
\end{pf}


\begin{cor} Suppose that $X$ is a locally compact left
$G$-space, and that $H$ is a closed normal subgroup of $G$ which
acts
freely and properly on $X$.  Then there is an isomorphism of
$\br_{G/H}(\hmx)$ onto
$\br_G(X)$ taking
$[B,\beta]$ to $[p^*(B),p^*(\beta)]=[p^*(B),\tau\tensor\beta]$.
\end{cor}
\begin{pf}
View $Y=X\times G/H$ as a left $G$-space via the diagonal action, and a
right $G/H$-space via right translation on the second factor.
Both actions are free, and the second action is proper.
To see that the first action is proper, suppose that $(x_n,t_nH)\to
(x,tH)$ while $(s_n\cdot x_n,s_nt_nH)\to (y,rH)$.
Then $s_nH\to sH$ for some $s\in G$.  Passing to a subsequence and
relabeling, we can assume that there are $h_n\in H$ such that $h_n
s_n\to s$ in $G$
.  But then $s_n \cdot
x_n\to y$ while $h_n\cdot(s_n\cdot x_n)\to s\cdot x$.
Since the $H$-action is proper, we can assume that $h_n\to h$ in $H$.
Thus $s_n\to h^{-1}s$, and this proves the claim.

The map $G\cdot(x,tH)\mapsto Ht^{-1}\cdot x$ is a
bijection $\phi$  of $G\backslash Y$ onto $H\backslash X$.
Further, $G\backslash Y$ is a right $G/H$-space and $H\backslash
X$ is a left $G/H$-space with
$$\phi\(v\cdot (s^{-1}H)\)=s H\cdot \phi(v).$$
(That is, $\phi$ is equivariant when the $G/H$-action on $G\backslash
Y$ is viewed as a left-action.)
Therefore,
\begin{equation}
\label{eq:*}
\br_{G/H}(G\backslash Y)\cong \br_{G/H}(H\backslash X)
\end{equation}
Similarly, $Y/(G/H)$ and $ X$ are isomorphic as left $G$-spaces
so that
\begin{equation}
\label{eq:**}
\br_G\(Y/(G/H)\)\cong \br_G(X).
\end{equation}
Finally, \thmref{thm:mt} implies that
\begin{equation}
\label{eq:dag}
\br_G\(Y/(G/H)\)\cong \br_{G/H}(G\backslash Y).
\end{equation}
Thus, Equations~(\ref{eq:*})--(\ref{eq:dag}) imply that there is an
isomorphism of $\br_{G/H}(H\backslash X)$
onto $\br_G(X)$ sending $(B,\beta)$ to
$\(\Ind_{G/H}^{X\times G/H}(B,\beta)/J,\tau\tensor\id\)$ with
$$J=\set{
f\in \Ind_{G/H}^{X\times G/H}(B,\beta)
:\text{$f(x,rH)(Hr^{-1}\cdot x)=0$ for all $x\in X$}}.$$
Define $\Phi:\Ind_{G/H}^{X\times G/H}(B,\beta)\to C_0(X,B)$ by $\Phi(
f)(x)=f(x,H)$.  Then $\Phi$ is onto (see, for example,
the first sentence of the proof of \cite[Lemma~2.6]{rw}).
Since
\begin{align*}
\Phi\(\tau_s\tensor\id(f)\)(x)&=\tau_s\tensor\id(f)(x,H)
=f(s^{-1}\cdot x, s^{-1} H) = \beta_{sH}\(f (s^{-1} \cdot x, H)\)
\\
&=\tau_s \tensor \beta_{sH}\(\Phi(f)\)(x),
\end{align*}
$\Phi$ is equivariant, and it only remains to show
that $\Phi$ induces a bijection of the quotient by $J$ with the
quotient of $C_0(X,B)$ by the balancing ideal~$I$.

However,
if $\Phi(f)\in I$, then $f(x,H)(H\cdot x)=0$ for all $x\in X$.
But then $f(x,rH)(Hr^{-1}\cdot x)=\beta_{rH}^{-1}\(f(x,H)\)(Hr^{-1}\cdot
x)$, which is zero since $\beta$ covers the $G/H$-action on $X$,
and $f\in J$.
The argument reverses, so $\Phi(J)= I$, and the result
follows.
\end{pf}

\def\mathcs{C^{\displaystyle *}} \def\cs{\ifmmode\mathcs\else$\mathcs$\fi}
\ifx\undefined\bysame
\newcommand{\bysame}{\leavevmode\hbox to3em{\hrulefill}\,}
\fi

\end{document}